\documentclass{svproc}

\usepackage{url}
\usepackage{epsfig,graphicx}
\usepackage{mathtools}
\usepackage{MnSymbol}

\newcommand{\loc}[1]{\zeta(#1)}
\DeclarePairedDelimiter{\ceil}{\lceil}{\rceil}

\begin{document}

\mainmatter

\title{Progress on the localization number of a graph}

\titlerunning{The localization number}  
%
\author{Anthony Bonato \and Melissa A.\ Huggan \and Trent G.\ Marbach}
\authorrunning{A. Bonato et al.} 
%
\tocauthor{Anthony Bonato,Melissa A.\ Huggan, Trent G.\ Marbach}
\institute{Ryerson University, Toronto, ON, Canada M5B2K3\\
\email{abonato@ryerson.ca,melissa.huggan@ryerson.ca,trent.marbach@gmail.com}}

\maketitle



\let\thefootnote\relax\footnotetext{The first author was supported by an NSERC Discovery Grant and the second author was supported by an NSERC Postdoctoral Fellowship.}

\begin{abstract}
We highlight new results on the localization number of a graph, a parameter derived from the localization graph searching game. After introducing the game and providing an overview of existing results, we describe recent results on the localization number.
We describe bounds or exact values of the localization number of incidence graphs of designs, polarity graphs, and Kneser graphs.
\end{abstract}

\section{Introduction to the localization game}
In graph searching, a set of pursuers attempts to locate, eliminate, or contain the threat posed by an evader in the network. The rules, specified from the outset, greatly determine the difficulty of the tasks proposed above. For example, the evader may be visible, but the pursuers may have limited movement speed, only moving to nearby vertices adjacent to them. Such a paradigm leads to the game of Cops and Robbers \cite{BN}, and deep questions like Meyniel's conjecture on the cop number of a graph.

We focus on a particular graph searching game in the present extended abstract. In the \emph{localization game}, there are two players moving on a reflexive graph, with one player controlling a set of $k$ \emph{cops}, where $k$ is a positive integer, and the second controlling a single \emph{robber}. Unlike in Cops and Robbers, the cops play with imperfect information: the robber is invisible to the cops during gameplay. The game is played over a sequence of discrete time-steps; a \emph{round} of the game is a move by the cops together with the subsequent move by the robber. The robber occupies a vertex of the graph, and when the robber is ready to move during a round, they may move to a neighboring vertex or remain on their current vertex. A move for the cops is a placement of cops on a set of vertices. Note that the cops are not limited to moving to neighboring vertices.

At the beginning of the game, the robber chooses their starting vertex. After this, the cops then make their move, followed by the robber; after that, the players move on alternate time-steps. Observe that any subset of cops may move in a given round. In each round, the cops occupy a set of vertices $u_1, u_2, \dots , u_k$ and each cop sends out a \emph{cop probe}, which gives their distance $d_i$, where $1\le i \le k$, from $u_i$ to the robber's vertex. The distances $d_i$ are either nonnegative integers or $\infty.$ Hence, in each round, the cops determine a \emph{distance vector} $(d_1, d_2, \dots , d_k)$ of cop probes, which is unique up to the ordering of the cops.
Given this distance vector, each vertex with distance $d_i$ to vertex $v_i$ for all $i$ is called a \emph{candidate}, as it may contain the robber. A candidate may not be unique.

\begin{figure}[ht!]
\begin{center}
\epsfig{figure=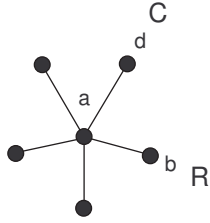,scale=1}
\caption{The localization game played on $K_{1,5}.$ The robber on $b$ cannot move to $a$ without being captured on the next cop turn. The cop moves to each leaf vertex and captures the robber in at most five moves.} \label{3lstar}
\end{center}
\end{figure}

The cops win if they have a strategy to determine, after a finite number of rounds, a unique candidate, at which time we say that the cops {\em capture} the robber. We also say the cops \emph{identify} or \emph{locate} the robber if there is a unique candidate. If there is no unique candidate in a given round, then the robber may move again. The cops may move to any other vertices in the next round resulting in an updated distance vector; meaning, the cops are not restricted to moving to adjacent vertices. The robber wins if they are never captured. See Figure~\ref{3lstar} for the game played on a star.

By randomly choosing vertices, the cops will eventually capture the robber. This strategy would let the cops win by chance, which is less interesting. We want to avoid this kind of nondeterministic capture, and assume that the robber is \emph{omniscient}, in the sense that they know the entire strategy for the cops. For example, with an omniscient robber who can anticipate the moves for the cop, we will need at least two cops to guarantee a capture in $K_3$. Otherwise, with one cop, there are always two candidates and the robber can escape capture.
From now on, we always assume we are playing against an omniscient robber.

For a graph $G$, define the \emph{localization number} of $G$, written $\loc{G}$, to be the least positive integer $k$ for which $k$ cops have a winning strategy. As placing a cop on each vertex gives a distance vector with unique value of $0$ on the location of the robber, the localization number is well-defined. If the order of $G$ is $n$, it is evident that $\loc{G} \le n-1.$

For disconnected graphs, unlike in Cops and Robbers, the cops are free to probe in different components. Eventually, they will locate the robber by probing a finite distance. Hence, the localization number of a disconnected graph is the maximum of its value on a component (rather than the sum, as is the case of the cop number). For this reason, we restrict our attention to connected graphs.

The localization game was first introduced for one cop by Seager~\cite{seager1}. In that paper, the robber was not allowed to \emph{backtrack}: they could not move to a vertex containing a cop in the previous round. Note that forbidding backtracking makes the game harder for the robber. The game in the present form was first considered by \cite{car}, who studied its properties on subdivisions of graphs. Localization was subsequently studied in several papers~\cite{BHM,BK,nisse1,BDELM,has,seager2}. We will highlight some of these results below.

The localization number is related to the metric dimension of a graph, in a way that is analogous to how the cop number is bounded above by the domination number. The \emph{metric dimension} of a graph $G$, written $\beta(G)$, is the minimum number of cops needed in the localization game so that the cops can win in one round; see \cite{hm,slater}. Hence, $\loc{G} \le \beta(G).$ However, there are cases where this inequality is far from tight.  For example, $\beta(K_{1,n}) = n-1,$ while $\loc{K_{1,n}} = 1.$

\section{Previous results on the localization number}\label{secpar3}

\begin{figure}[ht!]
\begin{center}
\epsfig{figure=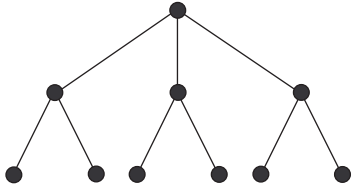,scale=1}
\caption{The tree $T_{3}.$} \label{t3}
\end{center}
\end{figure}

Trees were extensively studied by Seager~\cite{seager2} in early work on the localization number. For trees, we always have a bound on the localization number of 2.
\begin{theorem}[\cite{nisse2,seager2}]\label{3ltreet}
If $T$ is a tree, then $\loc{T} \le 2.$ We have that $\loc{T} = 2$ if and only if $T$ contains $T_3$ as a subgraph, as depicted in Figure~\ref{t3}.
\end{theorem}

In \cite{seager1}, a variant of the localization game was considered where the robber cannot visit a vertex visited by the cop in the previous round. We call this the no-backtrack condition and we refer to the corresponding parameter as $\zeta^*.$ Note that this makes the game harder for the robber, as they have fewer moves available to them. Assuming no-backtracking, it was shown in \cite{john} that if $\zeta^*(G)=1$, then $\chi(G) \le 4.$ There is a connection between the localization number and the chromatic number when backtracking is allowed.  In \cite{nisse1}, it was conjectured (see Conjecture 16 in the paper) that there is a positive-integer-valued function $f$ such that every graph $G$ with $\loc{G}  \le  k$ satisfies $\chi(G) \le  f (k).$ This conjecture was proved by Bonato and Kinnersley in \cite{BK} in the following theorem, by exploiting properties of graph degeneracy.
\begin{theorem}[\cite{BK}]\label{cor:chromatic}
For every graph $G$, we have $\chi(G) \le 3^{\loc{G}}$.
\end{theorem}

Besides the chromatic number, there are other upper bounds related to common graph parameters. We reference a nice connection between the localization number and the maximum degree of a graph.
\begin{theorem}[\cite{has}]\label{3lmax}
For a graph $G$ with maximum degree $\Delta,$ we have that
$$\loc{G}\le  \left \lfloor \frac{(\Delta +1)^2}{4} \right \rfloor + 1.$$
Further, if $\Delta =3,$ then $\loc{G} \le 3.$
\end{theorem}

There is also a bound given by pathwidth. We denote the \emph{pathwidth} of a graph $G$ by pw($G$).
\begin{theorem}[\cite{nisse1}]\label{3lpw1}
If $G$ is a graph, then $\loc{G} \le \mathrm{pw}(G).$
\end{theorem}

In \cite{BK}, the localization number of an $n$-dimensional hypercube, $Q_{n}$, was determined up to one of a small number of possible values.
\begin{theorem}[\cite{BK}]\label{3thm:hypercube}
For all positive integers $n$, we have that $$\ceil{\log_2 n} \le \loc{Q_n} \le \ceil{\log_2 n} + 2.$$
\end{theorem}

\section{The localization number of incidence graphs and graphs with diameter 2}\label{secfam3}

If $P$ is a projective plane of order $q,$ we use the notation $G(P)$ for its incidence graph.  We use the notation $X$ for the points of $P$ and $\mathcal{B}$ for the lines (or blocks) of $P$. We always assume that $q$ is a positive integer.

The \emph{robber territory} is defined as follows on a graph $G$. The robber territory is initialized to be $T_{0} = V(G)$. After the cops have moved on the cops $i$th turn to move, we define $T'_i$ to contain those vertices that are in $T_{i-1}$ or the neighbors of a vertex in $T_{i-1}$.  The vertices in $T'_i$ can be partitioned into classes such that each class contains exactly those vertices of $T'_i$ with identical distance vectors. The class that the robber currently resides on is defined as $T_i$.  As the robber has perfect information,  the robber is able to choose which of the classes of $T'_i$ is used for $T_i$.

For incidence graphs of projective planes, the localization number is known exactly. The proof uses a careful analysis of the robber territory, and how the cops may reduce it inductively.
\begin{theorem}[\cite{BHM}]\label{3lppp}
If $P$ is a projective plane of order $q,$ then $\loc{G(P)} = q+1.$
\end{theorem}
The localization number of other designs such as affine planes and Steiner triple systems were considered in \cite{BHM}. We summarize our results for the localization numbers of designs in the chart below. All the graphs in the chart are incidence graphs $G$ of designs, including balanced incomplete block designs (represented by \rm{BIBD}$(v,b,r,k,1)$) and Steiner systems (represented by \rm{STS}$(v)$). The columns list the design, bounds or exact values of $\zeta(G)$, and a reference to the appropriate theorem or corollary in \cite{BHM}.

\vspace{0.1in}
\begin{center}
\begin{tabular}{|l||c|c|}
\hline
\textbf{Design} & \textbf{Bounds or values} & \textbf{References in \cite{BHM}}\\
\hline\hline
\rm{BIBD}$(v,b,r,k,1)$ & $\zeta(G) \le 2r+k-3$ & Corollary~5 \\ \hline
Symmetric \rm{BIBD}$(v,b,r,k,1)$ & $\zeta(G) = k$ & Theorem~8 \\ \hline
Affine plane of order $q$ & $\zeta(G) =q$ & Theorem~9 \\ \hline
\rm{STS}$(v)$, $v>9$ & $\lfloor \frac{v-2}{8} \rfloor \le \zeta(G) \le \frac{v+1}{2}$ & Corollary~10, Theorem~11 \\ \hline
\rm{STS}$(v)$ & $\zeta(G) \le (1+o(1))v/3$ &  Theorem~12 \\ \hline
\end{tabular}
\end{center}

\vspace{0.1in}

There are other graphs we may define from projective planes that are distinct from incidence graphs. For a given projective plane of order $q$ with points $X$ and lines $\mathcal{B}$, a \emph{polarity} $\pi :X \rightarrow \mathcal{B}$ is a bijective mapping of points to lines such that $v\in \pi(u)$ whenever $u\in\pi(v)$.
The \emph{polarity graphs} are formed on vertex set $X$ by joining distinct vertices $u$ and $v$ if $u\in \pi (v)$ and $u \neq v$. Polarity graphs have $q^2+q+1$ vertices, and each vertex has degree $q$ or $q+1$. Such graphs are known to be $C_4$-free and diameter 2.

Polarity graphs were studied for the game of Cops and Robbers in~\cite{BB}, where bounds were given on their cop number. The following theorem was proven for the localization number of polarity graphs; the upper bound is the metric dimension.
\begin{theorem}[\cite{BHM1}]\label{finalt}
If $G$ is a polarity graph with order $q^2+q+1,$ then \[\frac{2q-5}{3} \leq \zeta(G) \leq 2q-1.\]
\end{theorem}

%
%

Kneser graphs are a well-known family of graphs. 
For integers $k,n \ge 1$ with $n > k,$ the \emph{Kneser graph} $K(k,n)$ has vertices labeled by the $k$-tuples on $\{ 1,2, \dots, n \},$ with two vertices adjacent if and only if their vertex labels are disjoint. We focus
on the diameter 2 case, which occurs exactly when $n \ge 3k.$ For each fixed even $k\geq 4$, the following theorem gives the exact value of the localization number up to an additive constant of an infinite subclass of Kneser graphs. The proof employs a new approach using hypergraph detection that bounds the localization number and metric dimension.
\begin{theorem}[\cite{BHM1}]\label{3finall}
For the localization number of Kneser graphs, we have the following bounds.
\begin{enumerate}
\item For a fixed even integer $k \geq 4$ and $n$ with $n \geq 3k$, we have that
\[\zeta(K(k,n)) =  \frac{n}{2} + \frac{n}{k} +O(1).\]
\item For a fixed odd integer $k \geq 3$ and $n$ with $n \geq 3k$, we have that
\[\frac{n}{2}+\frac{n}{k} -\frac{k}{2}-1 \leq \zeta(K(k,n)) \leq   \frac{n}{2} + \frac{n}{k}+ \frac{n}{2k} +O(1).\]
\end{enumerate}
\end{theorem}

We finish with open problems. An interesting question is to find tight bounds on the localization number of families of Steiner systems. The localization number has yet to be considered for block intersection graphs, point graphs, or Latin square graphs.  We do not know whether the bounds for the localization number of polarity graphs in Theorem~\ref{finalt} are sharp, and the same problem holds for the bounds in Theorem~\ref{3finall}. It would also be interesting to determine the localization number of Kneser graphs with diameter at least 3.

\end{document}